\documentclass[10pt]{article}

\usepackage{latexsym}
\usepackage{amsfonts}
\usepackage{amsmath}
\usepackage{amssymb}
\usepackage{textcomp}
\newtheorem{thm}{Theorem}
\newtheorem{lem}[thm]{Lemma}

\def\eps{\varepsilon}

\def\qed{\hfill $\vcenter{\hrule height .3mm
\hbox {\vrule width .3mm height 2.1mm \kern 2mm \vrule width .3mm
height 2.1mm} \hrule height .3mm}$ \bigskip}
\def\P{\mathbb{P}}

\def\EE{\mathbb{E}}

\def\RR{\mathbb{R}}

\def\Sph{S^{n-1}}

\title{\textsf{A Polynomial Number of Random Points does not Determine the Volume of a Convex Body} }
\author{Ronen Eldan \thanks{Partially supported by the Israel Science Foundation and by a Farajun Foundation Fellowship}}
\date{}
\begin{document}
\maketitle

\begin{abstract}
We show that there is no algorithm which, provided a polynomial number of random points uniformly distributed over a convex body
in $\RR^n$, can approximate the volume of the body up to a constant factor with high probability.
\end{abstract}

\section{Introduction}

Volume-related properties of high-dimensional convex bodies is one of the main topics of convex geometry in research today.
Naturally, calculating or approximating the volume of a convex body is an important problem.
Starting from the 1980's, several works have been made in the area of finding a fast algorithm for computing the volume of a convex body (see for example \cite{B},\cite{BF},\cite{LS},\cite{DFK},\cite{LV} and references therein). \\

These algorithms usually assume that the convex body $K \subset \RR^n$, is given by a certain oracle. An oracle is a "black box" which provides the
algorithm some information about the body. One example of an oracle is the \textbf{membership
oracle}, which, given a point $x \in \RR^n$, answers either "$x \in K$" or "$x \notin K$". Another example, is the \textbf{random point oracle}, which generates random points uniformly distributed over $K$. \\

All volume computing algorithms, known to the author, which appear in the literature use the membership oracle.
This note deals with a question asked by L. Lov\'{a}sz about the random point oracle.
It has been an open problem for a while whether or not it is possible to find a fast algorithm which
computes the volume of $K$ provided access to the random point oracle (\cite{GR}, \cite{L0}). \\

We answer this question negatively. In order to formulate our main result, we begin with some definitions. \\
An algorithm which uses the random point oracle is a (possibly randomized) function whose input is a finite sequence of random points generated according to the uniform measure
on $K$ and whose output is number, which is presumed be an approximation for the volume of $K$. The complexity of the algorithm will be defined by
the length of the sequence of random points. We are interested in the existence of algorithms with a complexity which depends polynomially on the dimension $n$. \\
We say that an algorithm is correct up to $C$ with probability $p$, if for any $K \subset \RR^n$, given the sequence of random points from
$K$, the output of the algorithm is between $\frac{Vol(K)}{C}$ and $C Vol(K)$, with probability at least $p$. \\
We prove the following theorem:
\begin{thm} \label{maingen}
There do not exist constants $C,p,\kappa>0$ such that for any dimension $n$, there exists an algorithm with complexity $O(n^\kappa)$ which is correct in
estimating the volume of convex bodies in $\RR^n$ up to $C$ with probability $p$.
\end{thm}

It is important to emphasize that this result is not a result in complexity theory. In this note we show that a polynomial number of points actually does not contain enough information to estimate the volume, regardless of the number of calculations, and hence, it is of information-theoretical nature. \\
For convex geometers, the main point in this study may be the additional information on volume distribution in convex
bodies it provides.  We suggest the reader to look this result in view of the recent results concerning the distribution of mass in convex bodies. In particular, results regarding thin-shell concentration and the Central Limit Theorem for Convex bodies, proved in
the general case by B. Klartag, show that essentially all of the mass of an isotropic convex body $K$ is
contained in a very thin-shell around the origin, and that almost all of the marginals are approximately gaussian. This may suggest that,
in some way, all convex bodies, when neglecting a small portion of the mass, behave more or less the same as a Euclidean ball in many senses.
Philosophically, one can also interpret these results as follows: provided a small number of points from a logarithmically-concave measure,
one cannot distinguish it from a spherically symmetric measure. For definitions and results see \cite{K}.
One of the main stages of our proof is to show that one cannot distinguish between
the uniform distribution over certain convex bodies, which are geometrically far from a Euclidean ball,
and some spherically symmetric distribution, when the number of sample points is at most polynomially large. \\

Here is a more quantitative formulation of what we prove:
\begin{thm} \label{main}
There exists a constant $\eps > 0$ and a number $N \in \mathbb{N}$ such that for all $n > N$, there does not exist an algorithm whose input is a sequence of length
$e^{n^\eps}$ of points generated randomly according to the uniform measure in a convex body $K \subset \mathbb{R}^n$, which determines $Vol(K)$ up to
$e^{n^{\eps}}$ with probability more than $e^{-n^\eps}$ to be correct.
\end{thm}

\textbf{Remark.} After showing that the volume of a convex body cannot be approximated, one may further ask: what about an algorithm that
estimates the \textbf{volume radius} of a convex body, defined by $VolRad(K) = Vol(K)^{\frac{1}{n}}$? A proof which shows that it is also
impossible has to be far more delicate than our proof. For example, under the hyperplane conjecture, it is easy
to estimate the volume radius of a convex body up to some $C>0$. \\

One may also compare this result to the two following related results: in a recent result N.Goyal and L.Rademacher (\cite{GR}) show that in order to \textbf{learn} a convex body, one needs at least $2^{c \sqrt \frac{n}{\eps}}$ random points. Learning a convex body rougly means finding a set having at most $\eps$ relative symmetric difference with the actual body (see \cite{GR}). Klivans, O'Donnel and Servedio (\cite{KOS}), show that any convex body can be agnostically learned with respect to the gaussian distribution using $2^{O(\sqrt n)}$ labelled gaussian samples. \\

The general idea of the proof is as follows.
Let $\{ K_\alpha \}_{\alpha \in I_1 }$ and $\{K_\alpha \}_{\alpha \in I_2}$ be two families of convex bodies.
For $i=1,2$, a probability measure $\mu_i$ on the set of indices $I_i$ induces a random construction of convex bodies, which in turn induces a probability measure $P_i$ on the set of sequences of points in $\RR^n$ in the following simple way: first generate an index $\alpha$ according to $\mu_i$, and then generate a sequence of $N$ uniformly distributed random
samples from $K_{\alpha}$. \\ 
In the proof we will define two distinct random constructions of convex bodies, $K_i = (\{ K_\alpha \}_{\alpha \in I_i }, \mu_i), ~~ i=1,2$ such that: \\
1. For every $\alpha_1 \in I_1$ and $\alpha_2 \in I_2$, the ratio between $Vol(K_{\alpha_1})$ and $Vol(K_{\alpha_2})$ is large. \\
2. If $N$ is not too large, both distributions $P_1, P_2$ are close in total variation distance to some distributions of samples in which the samples are independent and have a spherically symmetric law. \\
3. The radial profiles (hence the distribution of the Euclidean norm of a random sample) of typical random bodies $K_1,K_2$ are very close to each other. \\ \\
In other words, we will define two constructions of random convex bodies for which: 1. The typical volumes of the bodies they produce will be far from equal. 2. They will be both indistiguishable from spherically symmetric constructions for a polynomial number of samples. 3. The radial profiles they produce are indistiguishable from each other for a polynomial number of samples.  \\

To go on with the proof, a simple application of Yao's lemma will help us assume that the algorithm is deterministic.
A deterministic algorithm is actually a function $F : \mathbb{R}^{n^{\kappa + 1}} \to \mathbb{R}$ which takes a sequence
of points and returns the volume of the body. If the total variation distance between the
probabilities $P_1$ and $P_2$ defined above is small, then, there exists a set $A \subset \mathbb{R}^{n^{\kappa + 1}}$ which has a high probability with respect to
both $P_1$ and $P_2$. Obviously, for all $x \in A$, $F(x)$ is wrong in approximating the volume of at least one of the families. \\

In section 2, we will describe how we build these families of bodies, $\{K_{\alpha} \}$, using a random construction
which starts from a Euclidean ball, to which deletions which cut out parts of it, generated by some Poisson process, are applied.
Then, using elementary properties of the Poisson process and some concentration of measure properties of the ball,
we will see that the correlation between different points in polynomially long sequence of random points generated uniformly
from the body will be very weak (with respect to the generation of the body itself). Using this fact, we will only have
to inspect the distribution of a single random point. The construction will have a spherically-symmetric nature, so the density of a single random point will only depend
on its distance from the origin, and therefore we will only have to care about the distribution
of the distance of a point from the origin in the generated bodies.
The role of the following section, which is more technical but fairly delicate, will be to calibrate this construction so that these families have
different volumes, yet, approximately the same distribution of distance from the origin. \\ \\

Before we proceed to the proof, let us introduce some notation. In this note the number $n$ will always denote a dimension.
For an expression $f(n)$ which depends on $n$, by $f(n) = \mathbf{SE}(n)$
we mean: there exists some $n_0 \in \mathbb{N}$ and $\epsilon > 0$ such that for all $n > n_0$, $|f(n)| < e^{-n^{\epsilon}}$.
Also write $f(n) = g(n) (1 + \mathbf{SE}(n))$ for $\left |\frac{f(n)}{g(n)} - 1 \right | = \mathbf{SE}(n)$ and
$f(n) = g(n) + \mathbf{SE}(n)$ for $|f(n) - g(n)| = \mathbf{SE}(n)$.
The notation $f(n) \lesssim g(n)$ and $f(n) \gtrsim g(n)$ will be interpreted as $f(n) < g(n)$ and $f(n) > g(n)$ for $n$ large enough.

Moreover, we decide that $N=N(n)$, denotes the length of the sequence of random points. All throughout this note we assume that
there exists a universal constant $\eps > 0$, such that $N(n) < e^{n^\eps}$. \\ \\

\textbf{Acknowledgements} I am deeply grateful to my supervisor Prof. Bo\"{a}z Klartag for very useful discussions and encouragement all along my work on the subject.
I would also like to express my thanks to my supervisor, Prof Vitali Milman for introducing me to this question and encouraging me to work on it, to the referee for providing useful comments and insights, and to Prof Oded Regev for reviewing a preliminary version of the paper and providing numerous insightful comments and remarks. 

\section{The Deletion Process}
\label{sec2}
In this section we will describe the construction of the random bodies which will later be used as counter-examples. Our goal, after describing
the actual construction, will be to prove, using some simple properties of the Poisson distribution, a weak-correlation property between different
points generated from the body. \\
Denote by $D_n$ the $n$ dimensional Euclidean ball of unit radius, centered at the origin, and by $\omega_n$ its Lebesgue measure. \\
Recall that for two probability measures $P_1, P_2$ on a set $\Omega$, the total variation distance between the two measures is defined by
$$ d_{TV}(P_1,P_2) = \sup_{A \subseteq \Omega} |P_1(A) - P_2(A)| $$
One can easily check that if these measures are absolutely continuous with respect to some third measure $Q$, then it is also equal half the $L_1(Q)$
distance between the two densities. \\
Define $r_0 = n^{-\frac{1}{3}}$, and
$$
T_0(\theta) = D_n \cap \{ x; \langle x, \theta \rangle \leq r_0 \}.
$$

Let $T$ be a function from the unit sphere to the set of convex bodies, such that for every $\theta \in \Sph$,
$T(\theta)$ satisfies $T_0(\theta) \subseteq T(\theta) \subseteq D_n$.
(Recall that most of the mass of the Euclidean ball is contained in
$\{ x_1 \in [-1, C n^{-\frac 1 2}] \}$. So $T(\theta)$ contains almost all the mass of the Euclidean ball).
Moreover let $m>0$. We will now describe our construction of a random convex body, $K_{T, m}$. First, suppose that $m \in \mathbb{N}$.
Let $\Theta = (\theta_1, \theta_2, ..., \theta_m)$ be $m$ independent random directions, distributed according to the uniform measure on $\Sph$.
We define $K_{T,m}$ as,
$$ K_{T,m} = D_n \bigcap_i T(\theta_i). $$
Finally, instead of taking a fixed $m \in \mathbb{N}$, we take $\zeta$ to be a a Poisson random variable with expectation
$m$, independent of the above. We can now define $ K_{T, \zeta} $ in the same manner. \\
Let us denote the probability measure on the set of convex bodies induced by the process described above by $\mu$.
After generating the body $K_{T, m}$, which, from now on will be denoted just by $K$ wherever there is no confusion caused,
we consider the following probability space: let $\Omega = (D_n)^N$ be the set of sequences of length $N$ of points from $D_n$.
Denote by $\lambda$ the uniform probability measure on $\Omega$, and
for a convex body $K$ denote by $\lambda_K$ the uniform probability measure on $K^N = \prod_{1 \leq i \leq N} K \subseteq \Omega$.
Finally, define a probability measure $P=P_{T,m}$ on $\Omega$ as follows: for $A \subseteq \Omega$,
$$
P(A) = \int \lambda_K(A) d \mu(K) = \int \frac{Vol(K^N \cap A)}{Vol(K^N)} d \mu(K)
$$
(The measure $P$ describes the following process: first generate the random set $K$ according to construction described above, and then generate $N$ i.i.d random points, independent of the above, according
to the uniform measure on $K$).
Moreover, for $p = (x_1, ..., x_N) \in \Omega$, define $\pi_i(p) = x_i$, the projections onto the $i$-th copy of the Euclidean ball. \\
It it easy to check that $P$ is absolutely continuous with respect to $\lambda$. We define the following function on $\Omega$:
\begin{equation} \label{deff}
f_{T,m}(p) = \P(p \in K_{T,m}^N) = \P(\forall 1 \leq i \leq N, \pi_i(p) \in K_{T, m}).
\end{equation}
As we will see later, the function $f$ is related in a simple way to $\frac{dP}{d \lambda}$. Namely, we will have,
$$
\frac{dP}{d \lambda}(p) = (1 + \mathbf{SE}(n))  \frac{f(p)}{\int_{\Omega} f}
$$
for all $p$ in some subset of $\Omega$ with measure close to 1. For convenience, from now on $f_{T, m}$ will be denoted by $f$. \\ \\

We start with some simple geometric observations regarding $\Omega$. 
Denote by $\sigma$ the rotation invariant probability measure on
$\Sph$. Define, for $p \in \Omega$, $1 \leq i \leq N$,
\begin{equation} \label{defai}
A_i(p) = \{ \theta \in \Sph; \pi_i(p) \notin T(\theta) \}
\end{equation}
For $1 \leq i,j \leq N$, let $F_{i,j} \subset \Omega_N$ be the event, defined by
\begin{equation}
F_{i,j} = \left \{p; ~~ \frac{\sigma(A_{i}(p) \cap A_{j}(p))}{\sigma(A_{i}(p))} < e^{-n^{0.1}} \right \}
\end{equation}
and let,
\begin{equation}
F = \bigcap_{1 \leq i \neq j \leq N} F_{i,j}
\end{equation}
(which should be understood as "no two points are too close to each other", and, as we will see, will imply that points are
weakly correlated). We start with the following simple lemma.
\begin{lem} \label{mutualcaps} Under the above notations: \\
(i) $\lambda(F) = 1 + \mathbf{SE}(n) $. \\
(ii) There exists some $\eps_0 > 0$ such that:
if we assume that following condition holds,
\begin{equation} \label{volume}
\P_\mu (Vol(K) < \omega_n e^{-n^{\eps_0}}) < e^{-n^{\eps_0}}
\end{equation}
(hence, the volume of $K$ is typically not much smaller than the volume of $D_n$).
Then we have $P(F) = 1 + \mathbf{SE}(n)$. \\
\end{lem}
\textbf{Proof:} \\
(i) Let $p$ be uniformly distributed in $\Omega$. Denote $x_i = \pi_i(p)$, so
$x_1, x_2$ are independent points uniformly distributed in $D_n$. Let us calculate $\lambda(F_{1,2})$. \\
First, for a fixed $\theta \in \Sph$, one has
$$
\P(x_1 \notin T(\theta)) \leq \P(x_1 \notin T_0(\theta)) = \P(\{\langle x_1, \theta \rangle \geq r_0\})
$$
Recalling that $r_0 = n^{-\frac{1}{3}} \gg n^{-\frac{1}{2}}$,
by elementary calculations regarding marginals of the Euclidean ball, one gets
$$
\P(x_1 \notin T(\theta)) \lesssim e^{-n^{0.2}}
$$
Now, fix $x_2' \in D_n$. Define $A_i:=A_i(p)$. One has,
$$
\EE (\sigma (A_1 \cap A_2) | x_2 = x_2') = \int_{A_2} \P({\theta \in A_1}) d \sigma(\theta) = \int_{A_2} \P(x_1 \notin T(\theta)) d \sigma(\theta) \lesssim \sigma(A_2) e^{-n^{0.2}}
$$
And so,
\begin{equation}
\EE (\frac{\sigma (A_1 \cap A_2)}{\sigma(A_2)} | x_2 = x_2') \lesssim e^{-n^{0.2}}
\end{equation}
Now, this is true for every choice of $x_2'$, so integrating over $x_2'$ gives
$$
\EE \frac{\sigma (A_1 \cap A_2)}{\sigma(A_2)} \lesssim e^{-n^{0.2}}
$$
Now we use Markov's inequality to get
\begin{equation} \label{lambdaf}
\lambda(F_{1,2}^C) = \lambda(\left \{ \frac{\sigma(A_1\cap A_2)}{\sigma(A_2)} > e^{-n^{0.1}} \right \}) = \mathbf{SE}(n)
\end{equation}
A union bound completes the proof of (i). \\
\textbf{Proof of (ii)}
First, we can condition on the event $\{ Vol(K) > \omega_n e^{\eps_0} \}$ (with $\eps_0$ to be chosen later). (\ref{volume}) ensures us that it will happen with
probability $= 1 - \mathbf{SE}(n)$.
Observe that for any event $E \subset \Omega$ which is measurable by the $\sigma$-field generated by $\pi_1, \pi_2$, we have
\begin{equation} \label{twocoords}
\lambda_K (E) = \frac{\omega_n^2 \lambda((K \times K \times D_n \times ... \times D_n) \cap E)}{Vol(K)^2} \leq \frac{\omega_n^2 \lambda(E)}{Vol(K)^2}
\end{equation}
Now, taking $E=F_{1,2}^C$, choosing $\eps_0$ to be small enough and
using (\ref{lambdaf}) and (\ref{twocoords}) along with (\ref{volume}), one gets
$$
P(F_{1,2}) = 1 + \mathbf{SE}(n).
$$
Applying a union bound finishes the proof.
\qed \\ \\
We can now turn to the lemma which contains the main ideas of this section:
\begin{lem}: \label {indep}
There exist $\eps_0, \eps_1 > 0$ and $n_0$ such that for every $n > n_0$, the following holds:
Whenever $m$ is small enough such that the following condition is satisfied:
\begin{equation} \label{vol2}
\P (\{\theta \in K \}) > e^{-n^{\eps_0}}, ~~\forall \theta \in \Sph
\end{equation}
(hence, we are not removing too much volume, in expectation, even from the outer shell).
Then: \\
(i) We have,
\begin{equation} \label{volconc}
P(|Vol(K) - \EE(Vol(K))| > e^{-n^{\eps_1}} \EE(Vol(K))) = \mathbf{SE}(n)
\end{equation}
and also (\ref{volume}) holds. \\
(ii) For all $p \in F$, we have
$$ f(p) = (1 + \mathbf{SE}(n))  \prod_{j=1}^N \P(\pi_j(p) \in K)$$
In other words, if we define $\tilde f : D_n \to \RR$ as,
\begin{equation} \label{deftildef}
\tilde f (x) = \P(x \in K)
\end{equation}
then
\begin{equation} \label{fftilde}
 f(p) = (1 + \mathbf{SE}(n))  \prod_i \tilde f(\pi_i(p)), \forall p \in F.
\end{equation}
and,
$$
\mbox{(iii)}  ~~~ \frac{\EE(Vol(K^N \cap F))}{(\EE Vol(K))^N} - 1 = \mathbf{SE}(n)
$$
\end{lem}
\textbf{Proof}:
We begin by proving (ii). \\
Fix $p \in F$. Define $x_i = \pi_i(p)$, and $A_i = A_i(p) \subset \Sph$ as in (\ref{defai}). Also define
$G_j = \bigcap_{i \leq j} \{x_i \in K \}$.
Fix $2 \leq j \leq N$. Let us try to estimate $ P(G_j | G_{j-1})$. \\
When conditioning on the event $G_{j-1}$, we can consider our Poisson process as a superposition of three "disjoint" Poisson processes: the
first one, with intensity $\lambda_s$, only generates deletions that cut $x_j$, but leave all the $x_i$'s
for $i < j$ intact. The second one, with intensity $\lambda_u$ deletes $x_j$ along with one of
the other $x_i$'s, and the third one is the complement (hence, deletions that do not affect $x_j$). We have, recalling that the the expectation
of the number of deletions is $m$,
\begin{equation}
\lambda_s(\Sph) + \lambda_u(\Sph) =  m \sigma(A_j)
\end{equation}
Moreover,
\begin{equation}
\lambda_u(\Sph) \leq m \sum_{i < j} \sigma(A_i \cap A_j)
\end{equation}
(in the above formula we are including, multiple times, deletions that cut more than two points,
hence the inequality rather than equality). \\
Now, using the definition of $F$ one gets
\begin{equation}
\frac{\lambda_u(\Sph)}{\lambda_s(\Sph) + \lambda_u(\Sph)} = \mathbf{SE}(n)
\end{equation}
Note that (\ref{vol2}) implies
\begin{equation}
e^{-(\lambda_s(\Sph) + \lambda_u(\Sph))} \geq e^{-m \sigma(\{\theta; \frac{x_j}{|x_j|} \notin T(\theta)\}  )} \geq e^{-n^{\eps_0}}
\end{equation}
(the first inequality follows from the fact that $T(\theta)$ are star-shaped). The last two inequalities give,
\begin{equation}
\lambda_u(\Sph) = \mathbf{SE}(n)
\end{equation}
It follows that,
\begin{equation} \label{telescope1}
\left | \frac{P(G_j | G_{j-1})}{P(\{x_j \in K \})} - 1 \right | = \frac{e^{-\lambda_s(\Sph)}}{e^{-(\lambda_s(\Sph) + \lambda_u(\Sph))}} - 1 = \mathbf{SE}(n)
\end{equation}
Moreover, one has
\begin{equation} \label{telescope2}
P(G_N) = \prod_{j} P(G_j | G_{j-1})  = \prod_{j} \left (\frac{P(G_j | G_{j-1})}{P(\{x_j \in K \})} P(\{x_j \in K \}) \right)
\end{equation}
Using (\ref{telescope1}) and (\ref{telescope2}) we get
\begin{equation} \label{finalindep}
f(p) = P(G_N) = (1 + \mathbf{SE}(n))  \prod_j P(\{x_j \in K \})
\end{equation}
This proves (ii). \\
\textbf{Proof of (i):}
Showing that (\ref{volume}) holds is just a matter of noticing that $\P(x \in K)$ is monotone decreasing
with respect to $|x|$ and taking $\eps_0$ to be small enough. We turn to estimate $\EE(Vol(K)^2)$.
We have
\begin{equation}
\EE(Vol(K)^2) = \int_{D_n \times D_n} \P(\{x_1 \in K\} \cap \{x_2 \in K\}) dx1 dx2 =
\end{equation}
\begin{equation}
\int_{(D_n \times D_n) \cap F_{1,2}} \P(\{x_1 \in K\} \cap \{x_2 \in K\}) dx1 dx2 +
\end{equation}
$$
\int_{(D_n \times D_n) \cap F_{1,2}^C} \P(\{x_1 \in K\} \cap \{x_2 \in K\}) dx1 dx2
$$
(we will later see that the second summand is negligible).
Now, (\ref{finalindep}) gives
\begin{equation}
\int_{(D_n \times D_n) \cap F_{1,2}} \P(\{x_1 \in K\} \cap \{x_2 \in K\}) dx1 dx2 = 
\end{equation}
$$
(1 + \mathbf{SE}(n))  \int_{(D_n \times D_n) \cap F_{1,2}} \P(\{x_1 \in K\}) \P(\{x_2 \in K\}) dx1 dx2,
$$
which also implies that
$$\int_{(D_n \times D_n) \cap F_{1,2}} \P(\{x_1 \in K\} \cap \{x_2 \in K\}) dx1 dx2 > \frac{1}{2} e^{-2n^{\eps_0}} $$
Recall that $\lambda(F_{1,2}^C) = \mathbf{SE}(n)$ (as a result of the previous lemma). Taking $\eps_0$ to be small enough, we will get
$$
\EE(Vol(K)^2) = (1 + \mathbf{SE}(n))  \int_{(D_n \times D_n) \cap F_{1,2}} \P(\{x_1 \in K\} \cap \{x_2 \in K\}) dx1 dx2 =
$$
$$
(1 + \mathbf{SE}(n))  \int_{(D_n \times D_n) \cap F_{1,2}} \P(\{x_1 \in K\}) \P(\{x_2 \in K\}) dx1 dx2.
$$
On the other hand,
\begin{equation}
\EE (Vol(K))^2 = \int_{(D_n \times D_n)} \P(\{x_1 \in K\}) \P(\{x_2 \in K\}) dx1 dx2.
\end{equation}
Using the same considerations as above, the part of the integral over $F_{1,2}^C$ can be ignored, hence,
\begin{equation}
\EE (Vol(K))^2 = (1 + \mathbf{SE}(n))  \int_{(D_n \times D_n) \cap F_{1,2}} \P(\{x_1 \in K\}) \P(\{x_2 \in K\}) dx1 dx2.
\end{equation}
So we finally get
\begin{equation}
\EE (Vol(K)^2) = (1 + \mathbf{SE}(n))  \EE(Vol(K))^2
\end{equation}
Recalling that we assume (\ref{vol2}), using Chebishev's inequality, this easily implies (i), which finishes (ii). \\
For the proof of (iii),
$$
\EE (Vol(K^N \cap F)) = \int_F \P(p \in K^N) dp = (1 + \mathbf{SE}(n))  \int_F \prod_i \P(\pi_i(p) \in K) \leq (\EE Vol(K))^N.
$$
\qed \\ \\
Consider the density $\frac{dP}{d \lambda}$. Our next goal is to find a connection between this density and the function $f$.
Let $A \subseteq F \subset \Omega$. Using the concentration properties of $Vol(K)$, we will prove the following,
\begin{equation} \label{fandp}
P(A) = \frac{\int_A f(p) dp}{(\int_{D_n} \tilde f(x))^N} + \mathbf{SE}(n).
\end{equation}
where $f, \tilde f$ are defined in equations (\ref{deff}) and (\ref{deftildef}). \\
We have,
\begin{equation} \label{eqng0}
P(A) = \EE_{\mu} \left ( \frac{Vol(K^N \cap A)}{Vol(K^N)} \right )= \EE_{\mu} \left (\frac{Vol(K^N \cap A)}{Vol(K)^N} \right).
\end{equation}
By Fubini,
\begin{equation} \label{fubg}
\EE_{\mu} Vol(K^N \cap A) = \int_A f(p) dp.
\end{equation}
Consider the event
$$
G := \left \{ \left | \frac{Vol(K)^N}{\EE(Vol(K))^N} - 1 \right | < e^{-n^{\frac{\eps_1}{2}}} \right \}
$$
(where $\eps_1$ is the constant from lemma \ref{indep}). We have by the definition of $G$,
\begin{equation} \label{eqng1}
\int_{G} \frac{Vol(K^N \cap A) }{Vol(K)^N} d \mu(K) = \frac{\int_{G} Vol(K^N \cap A) d \mu(K)}{\EE(Vol(K))^N} + \mathbf{SE}(n).
\end{equation}
It follows from part (i) of lemma \ref{indep} that,
$$
\mu(G) = \P(\left| (\frac{Vol(K)}{\EE(Vol(K))})^N - 1 \right | \leq e^{-n^{\frac{\eps_1}{2}}}) \geq
$$
$$
\P(\left| \frac{Vol(K)}{\EE(Vol(K))} - 1 \right | \leq 2 N e^{-n^{\frac{\eps_1}{2}}}) \geq \P(\left| \frac{Vol(K)}{\EE(Vol(K))} - 1 \right | \leq e^{-n^{\eps_1}}) = 1 + \mathbf{SE}(n).
$$
So $\mu(G) = 1 + \mathbf{SE}(n)$ which gives,
\begin{equation} \label{eqng2}
\int_{G^C} \frac{Vol(K^N \cap A)}{Vol(K)^N} d \mu (K) \leq \mu(G^C) = \mathbf{SE}(n)
\end{equation}
We will also need:
\begin{equation} \label{eqng3}
\frac{ \int_{G^C} Vol(K^N \cap A) d \mu(K)}{(\EE Vol(K))^N} = \mathbf{SE}(n)
\end{equation}
To prove this, first recall that $A \subseteq F$. This gives,
\begin{equation} \label{eqng4}
\frac{ \int_{G^C} Vol(K^N \cap A) d \mu(K)}{(\EE Vol(K))^N} \leq \frac{ \int_{G^C} Vol(K^N \cap F) d \mu (K)}{(\EE Vol(K))^N} =
\end{equation}
$$
\frac{ \EE_\mu Vol(K^N \cap F)}{(\EE Vol(K))^N} - \frac{ \int_{G} Vol(K^N \cap F) d \mu (K)}{(\EE Vol(K))^N}.
$$
Now,
$$
\int_G \frac{Vol(K^N \cap F)}{Vol(K^N)} d \mu(K) = \EE_\mu \frac{Vol(K^N \cap F)}{Vol(K^N)} + \mathbf{SE}(n) = P(F) + \mathbf{SE}(n) = 1 + \mathbf{SE}(n)
$$
so,
\begin{equation} \label{eqng5}
\frac{ \int_{G} Vol(K^N \cap F) d \mu (K)}{(\EE Vol(K))^N} = 1 + \mathbf{SE}(n)
\end{equation}
Using part (iii) of lemma \ref{indep} along with
(\ref{eqng4}) and (\ref{eqng5}) proves $(\ref{eqng3})$. \\
Plugging together (\ref{eqng0}), (\ref{eqng1}), (\ref{eqng2}) and (\ref{eqng3}) imply
\begin{equation} \label{fandpi}
P(A) = \EE_{\mu} \frac{Vol(K^N \cap A)}{Vol(K)^N} = \int_G \frac{Vol(K^N \cap A)}{Vol(K)^N} d \mu (K) + \mathbf{SE}(n)
\end{equation}
$$
= \frac{\int_G Vol(K^N \cap A) d \mu (K)}{\EE(Vol(K))^N} + \mathbf{SE}(n) = \frac{\EE_{\mu} Vol(K^N \cap A) }{\EE(Vol(K))^N} + \mathbf{SE}(n)
$$
Recall that, as a result of Fubini's theorem,
\begin{equation} \label{expvolf}
\EE_\mu (Vol(K)) = \int_{D_n} \tilde f(x) dx.
\end{equation}
Plugging (\ref{fandpi}), (\ref{expvolf}) and (\ref{fubg}) proves (\ref{fandp}).
We would now like to use the result of lemma \ref{indep}, to replace $f$ with $\tilde f$.
Let $A' \subseteq \Omega$. Define $A=A' \cap F$,
$$
P(A') = P(A) + P(A' \cap F^C).
$$
Part (ii) of lemma \ref{mutualcaps} with (\ref{fandp}) gives
$$
P(A') = P(A) + \mathbf{SE}(n) = \frac{\int_{A} f(p) dp}{(\int_{D_n} \tilde f(x)dx)^N} + \mathbf{SE}(n).
$$
We can now plug in (\ref{fftilde}) to get
$$
P(A') = \frac{\int_{A} \prod_i \tilde f(\pi_i(p)) d p}{(\int_{D_n} \tilde f(x))^N} + \mathbf{SE}(n).
$$
So, finally defining
$$
\frac{d \tilde P}{d p} = \frac{\mathbf{1}_{\{ p \in F \}} \prod_i \tilde f(\pi_i(p))}{(\int_{D_n} \tilde f(x))^N} = \mathbf{1}_{\{ p \in F \}} \prod_i \frac{\tilde f(\pi_i(p))}{\int_{D_n} \tilde f(x) dx}
$$
we have proved the following lemma:
\begin{lem} \label{ptilde}
Suppose that the condition (\ref{vol2}) from Lemma \ref{indep} holds. Then one has
$$
d_{TV} (P, \tilde P) = \mathbf{SE}(n)
$$
\end{lem}
Note that the measure $\tilde P$ is not, in general, a probability measure. The lemma, however,
ensures us that $\tilde P(\Omega)$ is very close to 1. \\
Recall that our plan is to find two families of convex bodies, which are achieved by two pairs $(T_1, m_1)$ and $(T_2, m_2)$,
such that $d_{TV} (P_1, P_2)$ is small, even though their volumes differ. \\
The above lemma motivates us to try to find such pairs with $\frac{\tilde f_1}{\int \tilde f_1} = \frac{\tilde f_2}{\int \tilde f_2} + \mathbf{SE}(n)$.
We formulate this accurately in the following lemma.

\begin{lem} \label{sec1final}
Suppose there exist two pairs $(T_i, m_i)$ for $i=1,2$ such that (\ref{vol2}) is satisfied, and in addition,
defining $\tilde f_1$ and $\tilde f_2$ as in (\ref{deftildef}),
\begin{equation} \label{l1dist}
\left | \left | \frac{\tilde f_1}{\int_{D_n} \tilde f_1} - \frac{\tilde f_2}{\int_{D_n} \tilde f_2} \right | \right |_{L_1(D_n)} = \mathbf{SE}(n)
\end{equation}
Then $d_{TV}(P_1, P_2) = \mathbf{SE}(n)$.
\end{lem}
\textbf{Proof}: \\
Using the previous lemma, it is enough to show that $d_{TV}(\tilde P_1, \tilde P_2) = \mathbf{SE}(n)$.
Define $g_i = \frac{\tilde f_i}{\int_{D_n} \tilde f_i}$. We have
$$
d_{TV} (\tilde P_1, \tilde P_2) \leq \int_{\Omega} \left | \prod_{1 \leq i \leq N} g_1(\pi_i(p)) - \prod_{1 \leq i \leq N} g_2(\pi_i(p)) \right | \leq
$$
$$
 \sum_{1 \leq j \leq N} \int_{\Omega} \left |\prod_{1 \leq i \leq j} g_1(\pi_i(p)) \prod_{j+1 \leq i \leq N} g_2(\pi_i(p)) - \prod_{1 \leq i \leq j+1} g_1(\pi_i(p)) \prod_{j+2 \leq i \leq N} g_2(\pi_i(p)) \right | =
$$
$$
N \int_{D_n} |g_1(x) - g_2(x)| = \mathbf{SE}(n)
$$
\qed

In the next section we deal with how to calibrate $T_i$ and $m_i$ so that (\ref{l1dist}) holds.

\section{Building the two profiles}

Our goal in this section is to build convex bodies with a prescribed radial profile. \\
For a measurable body $L \subset \RR^n$, define
\begin{equation} \label{tvoldist}
g_{L}(r) = 1 - \sigma(\frac{1}{r} L \cap \Sph),
\end{equation}
This function should be understood as the "profile" of mass of the complement of $L$, which will eventually be
 the ratio of mass which a single deletion removes, in expectation, as a function of the distance from the origin. Define $g_i(r) = g_{T_i}(r)$. \\

Let us try understand exactly what kind of construction we require. Fix $x \in D_n$. Keeping in mind that the function $T_i(\theta)$ commutes with orthogonal transformations, we learn that the probability that $x$ is removed in a single deletion of $T_i$ is exactly $g_i(|x|)$. By elementary properties of the Poisson process, this gives,
\begin{equation} \label{poissonvol}
\P(x \in K_i) = \exp[-m_i g_i(|x|)].
\end{equation} 
In view of (\ref{l1dist}), we would like the ratio $ \frac{\P(x \in K_1)}{\P(x \in K_2)}$ to be (approximately) constant. Using (\ref{poissonvol}), we see that the latter follows from
$$ m_1 g_1 (|x|) - m_2 g_2 (|x|) = C. $$ \\
If we choose to pick $m_2 = 2 m_1$, this equality will be implied by the following requirements on $T_1, T_2$:
\begin{equation} \label{derratio1}
g_1(1) = g_2(1) \neq 0, ~~ \mbox{ and } ~~ g_1'(r) = 2 g_2'(r), ~~ r \in[0,1].
\end{equation}
Assuming (\ref{volume}) holds and making use of the concentration of the radial profile of $D_n$, we will actually only be required to make sure the derivatives are proportional for $r \in [1 - n^{-0.99}, 1]$.   \\
Note that when (\ref{derratio1}) is attained, by picking different values of $m_1$, the ratio between the expected volumes of $K_1$ and $K_2$ can be made arbitrarily large while the expected radial profiles remain about as close. Lemma (\ref{sec1final}) will then ensure us that this is enough for the distributions to be indistiguishable. \\ \\
The above is established in the main lemma of this section:
\begin{lem} \label{bodies}
For every dimension $n$, there exist two convex bodies $T_1, T_2 \subset \RR^n$, satisfying the following:
\begin{equation} \label{tsupset}
\mbox{(i) } D_n \supseteq T_i \supseteq D_n \cap \{x; \langle x, e_1 \rangle \leq n^{-\frac{1}{3}} \}, ~~ i=1,2
\end{equation}
(ii) The radial profiles satisfy,
\begin{equation} \label{derratio}
g_1(1) = g_2(1) \neq 0, ~~ \mbox{ and } ~~ g_1'(r) = 2 g_2'(r) ~~ \forall r \in [1 - n^{-0.99}, 1]
\end{equation}
\end{lem}
\medskip
To achieve this, we begin by describing the following construction:
Define $\delta_0 = n^{-\frac{1}{4}}$, $\delta_1 = n^{-0.99}$. For every two constants $a, b$ such that $a \in [2,200]$ and $b \in [-1000,1000]$,
let $f=f_{a,b}$ be the linear function with negative slope which satisfies:
\begin{equation}\label{lin1}
f (\delta_0 (1 + \delta_1b)) = \sqrt{1 - (\delta_0 (1 + \delta_1b))^2}
\end{equation}
and,
\begin{equation} \label{lin2}
\min_{x \in \RR} \sqrt{x^2 + f^2(x)} = a \delta_0
\end{equation}
(hence, it is a line of distance $a \delta_0$ from the origin which meets the unit circle at $x=\delta_0(1 + b \delta_1)$. Note that there exists
such a linear function with negative slope since $a \delta_0 \gg \delta_0(1 + b \delta_1)$). We define a convex body $T_{a,b}$ by,
\begin{equation}
T_{a,b} = D_n \cap \left \{(x, \vec y) \in \RR \times \RR^{n-1} = \RR^n ; |y| \leq f(x) \right \}
\end{equation}
(an intersection of the ball with a cone defined by a linear equation the coefficients of which depend of $a,b$). \\
Recall that we require that $a > 2$ and $b > -1000$. First of all, it follows directly from requirement (\ref{lin1}) and from the fact that
the slope of $f$ is negative, that $T_{a,b}$ satisfies (\ref{tsupset}) (since $\delta_0 \gg n^{-1/3}$). \\
Define $g_{a,b}(r) = g_{T_{a,b}}(r)$ as in (\ref{tvoldist}). Let us find an expression for $g_{a,b}(r)$.
First, a simple calculation shows that (\ref{lin2}) implies that the function $f_{a,b}$ intersects the $x$ axis at $x < 2 a \delta_0$.
This shows that $T_{a,b} \cap r \Sph$ has only one connected component for all $r > \frac{1}{2}$ (hence, the vertex of the cone is inside the sphere). \\
Consider the intersection $\frac{1}{r} T_{a,b} \cap \Sph$. If $r > \frac{1}{2}$, it must be a set of the form $\Sph \cap \{x_1 < x(a,b,r)\}$, for some function $x(a,b,r)$.
Let us try to find the expression for this function.
Equation (\ref{lin2}) shows that $T_{a,b}$ is an intersection of $D_n$ with halfspaces at distance $a \delta_0$ from the origin. This implies that
$x(a,b,r)$ must satisfy
$$
x(a,b,r) = \sin(\arcsin(\frac{a \delta_0}{r}) + c)
$$
for some constant $c$ (draw a picture). To find the value of $c$, we use (\ref{lin1}) to get $x(a,b,1) = \delta_0 (1 + b \delta_1)$, and so
\begin{equation} \label{finalxabr}
x(a,b,r) = \sin(\arcsin(\frac{a \delta_0}{r}) - \arcsin(a \delta_0) + \arcsin(\delta_0 (1 + b \delta_1))).
\end{equation}
Next, define
$$ \Psi(x) = \frac{1}{\omega_n} \int_{\min(x,1)}^1 (1 - t^2)^{\frac{n-3}{2}} dt,$$
the surface area measure of a cap the base of which has distance $x$ from the origin.
We have finally,
\begin{equation}
g_{a,b}(r) = \sigma(\Sph \cap \{x_1 \geq x(a,b,r) \}) = \Psi(x(a,b,r)).
\end{equation}
Given a subset $I' \subseteq \RR \times \RR$, we define
\begin{equation} \label{constrk}
K_{I'} = \bigcap_{(a,b) \in I'} T_{a,b}.
\end{equation}
Clearly
$$
g_{I'}(r) := g_{K_{I'}}(r) = \sup_{(a,b) \in I'} g_{a,b}(r)
$$
Our goal is to choose such a subset so that (\ref{derratio}) is fulfilled. We will use the following elementary result:

\begin{lem} \label{createconv}
Let $c > 0$, and let $\{ f_{\alpha} \}_{\alpha \in I}$ be a family of twice-differentiable functions defined on $[x_1,x_2]$ such that for every triplet
$(x, y, y') \in [x_1,x_2] \times [y_1,y_2] \times [y_1',y_2']$, there exists $\alpha \in I$
such that
\begin{equation}
f_{\alpha}(x) = y, ~~ f_{\alpha}(x)' = y', ~~ f''(t) \leq c, \forall t \in [x_1,x_2]
\end{equation}
then for every twice differentiable function $g: [x_1,x_2] \to [y_1,y_2]$ with
\begin{equation} \label{eqg}
g'(x) \in [y_1',y_2'], ~~ g''(x) > c,
\end{equation}
there exists a subset $I' \subset I$ such that
\begin{equation}
g(x) = \sup_{\alpha \in I'} f_{\alpha}(x)
\end{equation}
\end{lem}

In view of the above lemma, we would like to show that by choosing appropriate values of $a,b$, one can attain functions $g_{a,b}$
which, for a fixed $r_0$, have prescribed values $g_{a,b}(r_0), g_{a,b}'(r_0)$, and a small enough second derivative.

Define $r(u) = 1 - \delta_1 u$. Note that substituting $r \to u$, almost all of the mass of the Euclidean ball is contained in
$u \in [0,1]$ (the thin shell of the Euclidean ball).
We now turn to prove the following lemma:

\begin{lem} \label{uglylemma}
Suppose that $(u, g_0, g_0')$ satisfy $0 \leq u \leq 1$,
$$
\Psi(\delta_0) - 100 \delta_0 \delta_1 \Psi'(\delta_0) \leq g_0 \leq \Psi(\delta_0) + 100 \delta_0 \delta_1 \Psi'(\delta_0),
$$
$$
10 \delta_0 \delta_1 \Psi'(\delta_0) \leq g_0' \leq 100 \delta_0 \delta_1 \Psi'(\delta_0).
$$
There exist constants $a \in [2, 200], b \in [-1000, 1000]$ such that $g_{a,b}(r(u)) = g_0$, $(g_{a,b}(r(u)))' = g_0'$ and
$g_{a,b}(r(t))'' \leq \delta_0 \delta_1 \Psi'(\delta_0), \forall 0 \leq t \leq 1$.
\end{lem}
\textbf{Proof:}
Throughout this proof we always assume $u \in [0,1], a \in [2,200]$ and $b \in [-1000, 1000]$. \\
Let us inspect the function $x(a,b,r))$ defined in (\ref{finalxabr}).
Differentiating it twice, while recalling that $a \delta_0 \ll \frac{1}{2}$, gives us the following fact:
there exists $C > 0$ independent of $n$, such that
$ |\frac{\partial^2}{\partial r^2} x(a,b,r)| < C $.
Consider $x(a,b,u) := x(a,b,r(u))$. One has,
\begin{equation}
x_{uu}(a,b,u) = O(\delta_1^2)
\end{equation}
(here and afterwards, by "$O$", we mean that the term is smaller than some universal constant times the expression inside the brackets, which is valid as long as $u,a,b$ attain values in the intervals defined above).
This implies that for all $u \in [0,1]$,
$$
x_{u}(a,b,u) = x_{u}(a,b,0) + O(\delta_1^2) =
$$
$$
a \delta_0 \delta_1 \sin'(\arcsin(\delta_0 (1 + b \delta_1)))(1 + O(\delta_0))  + O(\delta_1^2)=
$$
$$
a \delta_0 \delta_1 (1 + O(\delta_0))
$$
and so,
$$
x(a,b,u) = x(a,b,0) + a \delta_0 \delta_1 u (1 + O(\delta_0)) = \delta_0 + \delta_0 \delta_1(a u + b)(1 + O(\delta_0))
$$
Let us now define $w(a,b,u) = \frac{1}{\delta_1} (\frac{x(a,b,r(u))}{\delta_0} - 1)$.
So,
\begin{equation}
w(a,b,u) = (au + b)(1 + O(\delta_0))
\end{equation}
and,
\begin{equation} \label{dersw}
w_u(a,b,u) = a(1 + O(\delta_0)), ~~ w_{uu}(a,b,u) = O(\frac{\delta_1}{\delta_0}).
\end{equation}
Next, we consider $g_{a,b}(r(u)) = \Psi(x(a,b,u)) = \Psi(\delta_0(1 + \delta_1 w(a,b,u))$. We have,
\begin{equation}
\Psi(x(a,b,u)) = \Psi(\delta_0) + \delta_0 \delta_1 \Psi'(\delta_0) w(a,b,u) + \frac{\delta_0^2 \delta_1^2}{2} \Psi''(t)w(a,b,u)^2
\end{equation}
for some $t \in [\delta_0, x(a,b,u)]$. But, note that the following holds,
\begin{equation} \label{secondderiv}
(\log \Psi'(v))' = \frac{\Psi''(v)}{\Psi'(v)} = - \frac{2v\frac{n-3}{2} (1 - v^2)^{\frac{n-5}{2}}}{(1 - v^2)^{\frac{n-3}{2}}} = - \frac{ v (n-3)}{(1 - v^2)},
\end{equation}
and for all $v \in [\frac{\delta_0}{2}, 2 \delta_0]$,
$$
(\log \Psi'(v))' = O(n \delta_0).
$$
Integration of this inequality yields that for $t$ such that $t - \delta_0 = O(\delta_0 \delta_1)$, one has
$$
\log \Psi'(t) - \log \Psi'(\delta_0) = O(n \delta_0^2 \delta_1)
$$
or,
\begin{equation} \label{firstderivt}
\Psi'(t) = \Psi'(\delta_0) (1 + O(n \delta_0^2 \delta_1))
\end{equation}
Combining (\ref{secondderiv}) and (\ref{firstderivt}) gives
\begin{equation} \label{secondderfinal}
\delta_0^2 \delta_1^2 \Psi''(t) = O(\Psi'(\delta_0) \delta_1^2 n \delta_0^3) = o(\Psi'(\delta_0) \delta_0 \delta_1).
\end{equation}
This finally gives,
\begin{equation} \label{psifinal}
g_{a,b}(r(u)) = \Psi(x(a,b,u)) = \Psi(\delta_0) + (\delta_0 \delta_1 \Psi'(\delta_0) w(a,b,u)) (1 + o(1))
\end{equation}
$$
= \Psi(\delta_0) + \delta_0 \delta_1 \Psi'(\delta_0)(au+b)(1 + o(1))
$$
Next we try to estimate the derivative of $\Psi(x(a,b,u))$. We have,
\begin{equation}
\frac{\partial }{ \partial u} g_{a,b}(r(u)) = \frac{\partial}{\partial u} \Psi(x(a,b,u)) =
\end{equation}
$$
\Psi'(x(a,b,u)) x_u(a,b,u) = \Psi'(x(a,b,u)) \delta_0 \delta_1 w_u(a,b,u)
$$
And using (\ref{firstderivt}),
\begin{equation} \label{derfinal}
\frac{\partial}{\partial u} \Psi(x(a,b,u)) = \Psi'(\delta_0) (1 + o(1)) \delta_0 \delta_1 w_u(a,b,u) =
\end{equation}
$$
(a \delta_0 \delta_1 \Psi'(\delta_0))(1 + o(1)).
$$
Using the continuity and $\Psi$ and $x(a,b,u)$, we can now conclude the following:
for any fixed $b \in [-1000,1000]$ and $u \in [0,1]$, an inspection of equation (\ref{derfinal}) teaches us that
when $a$ varies in $[2,200]$, $\frac{\partial}{\partial u} \Psi(x(a,b,u))$ can attain all values in the range
$[3 \delta_0 \delta_1 \Psi'(\delta_0), 100 \delta_0 \delta_1 \Psi'(\delta_0)]$. An inspection of equation (\ref{psifinal}) shows that afterwards, by letting $b$ vary in $[-1000,1000]$,
$g_{a,b}(r(u))$ will attain all values in $[\Psi(\delta_0) - 100 \delta_0 \delta_1 \Psi'(\delta_0), \Psi(\delta_0) + 100 \delta_0 \delta_1 \Psi'(\delta_0)]$.
To estimate the second derivative, $g_{a,b}''$, we write
$$
\frac{\partial^2}{\partial u^2} \Psi(x(a,b,u)) = \Psi''(x(a,b,u)) \delta_0^2 \delta_1^2 w_u^2(a,b,u) + \delta_0 \delta_1 \Psi'(x(a,b,u)) w_{uu}(a,b,u) w_{u}(a,b,u)
$$
(using (\ref{dersw}) and (\ref{secondderfinal}))
$$
= o(\delta_0 \delta_1 \Psi'(\delta_0)) + O(\delta_1^2 \Psi'(\delta_0)) = o(\delta_0 \delta_1 \Psi'(\delta_0))
$$
This completes the lemma. \qed

We are now ready to prove the main lemma of the secion. \\
\textbf{Proof of lemma \ref{bodies}} \\
Define:
$$
f_i(r) = \Psi(\delta_0) + C_i \delta_0 \delta_1 \Psi'(\delta_0) (u + 1)^2
$$
with $C_1 = 20, C_2 = 40$. Usage of lemmas (\ref{uglylemma}) and (\ref{createconv}) shows that there exist two subsets $I_1, I_2$
of $[2,200] \times [-1000,1000]$ such that the bodies $T_i = T_{I_i}$ that we constructed in (\ref{constrk})
satisfy (\ref{derratio}). Also, (\ref{tsupset}) is satisfied, since it is satisfied for $T_{a,b}$ for all
$(a,b) \in [2,200] \times [-1000,1000]$, as we have seen.
\qed

\section{Tying up Loose Ends}
\textbf{Proof of theorem \ref{main}}: \\
Use lemma \ref{bodies} two build the two bodies $T_i$. Let $U_{\theta}$ be an orthogonal transformation which sends $e_1$ to $\theta$.
Define $T_i(\theta) = U_{\theta}(T_i)$ (note that the choice of orthogonal transformation does not matter because
$T_i$ are bodies of revolution around $e_1$). Define the functions $g_i = g_{T_i}$ as in (\ref{tvoldist}).
Let $m_1 = \frac{n^{\eps}}{g_1(1)}$, with $\eps > 0$ to be chosen later. Define $m_2 = 2 m_1$. So, (\ref{derratio}) implies that
\begin{equation} \label{volratio}
m_2 g_2(r) = m_1 g_1(r) + m_1 g_1(1), ~~ \forall r \in [1 - n^{-0.99}, 1].
\end{equation}
Now, let $K_i = K_{T_i, m_i}$ be the random bodies we constructed in section 2. \\
For a fixed $x \in D_n$, as in (\ref{poissonvol}), we have,
\begin{equation} \label{pxi}
\tilde f_i(x) = P(x \in K_i) = e^{-m_i \sigma(\{x \notin T_i(\theta)\})} = e^{- m_i g_i(|x|)}.
\end{equation}
Now, (\ref{volratio}) and (\ref{pxi}) give
\begin{equation} \label{ratioff}
\frac{\tilde f_1(x)}{\tilde f_2(x)} = e^{m_1 (g(1) )} = e^{n^{\eps}}
\end{equation}
for all $x$ with $|x| \in [1 - n^{-0.99}, 1]$. \\
Let us choose $\eps$ to be small enough so that
$$
m_2 g_2(1) < n^{\eps_0}
$$
where $\eps_0$ is the constant from (\ref{vol2}). Clearly, that ensures that (\ref{vol2}) holds for both random
bodies $K_i$. Now, $\eps$ can be made further smaller, so that concentration properties of the Euclidean ball will give us,
\begin{equation} \label{thinshellf}
\int_{D_n} \tilde f_i = (1 + \mathbf{SE}(n))  \int_{D_n \setminus (1 - n^{-0.99}) D_n} \tilde f_i
\end{equation}
for $i=1,2$. Clearly, the above can still be satisfied for some universal constant $\eps > 0$ as long as $n$ is large enough.
Next, (\ref{ratioff}) and (\ref{thinshellf}) imply that
$$
\frac{\int_{D_n} \tilde f_1}{\int_{D_n} \tilde f_2} = (1 + \mathbf{SE}(n))  e^{n^\eps}
$$
and so (again, taking $\eps$ to be small enough) one gets
$$
\int_{D_n} \left |\frac{\tilde f_1}{\int_{D_n} \tilde f_1} - \frac{\tilde f_2}{\int_{D_n} \tilde f_2} \right | dx = \int_{(1 - n^{-0.99}) D_n} \left |\frac{\tilde f_1}{\int_{D_n} \tilde f_1} - \frac{\tilde f_2}{\int_{D_n} \tilde f_2} \right | dx + \mathbf{SE}(n) = \mathbf{SE}(n).
$$
Now use Lemma \ref{sec1final} to get that
\begin{equation} \label{totalvarsmall}
d_{TV} (P_1, P_2) = \mathbf{SE}(n)
\end{equation}
Denote $R = \frac{1}{2} e^{n^{\eps}}$. Then,
$$
\EE(Vol(K_1)) = (1 + \mathbf{SE}(n))  2 R \EE(Vol(K_2)).
$$
Suppose by negation that there exists a classification function $F: \omega \to \mathbb{R}$ that determines the volume of a body $K$ up to a constant
$e^{n^{\eps_2}}$ with probability $0.52$. Denote $L = [\frac{\EE(Vol(K_1))}{R}, R \EE(Vol(K_1))]$.
Note that using (\ref{volconc}), the "correctness" of the function implies that
$$
P_1(F(p) \in L) \geq 0.51
$$
Denote $A \subset \Omega$ as $A = \{p \in \Omega; F(p) \in L\}$. Then $P_1(A) > 0.51$, and
(\ref{totalvarsmall}) implies that also $P_2(A) > 0.51$. But this means that
$$
P_2(F(p) \in L) > 0.5
$$
But clearly, again, (\ref{volconc}) implies that with probability $= 1 + \mathbf{SE}(n)$, the volume of $K_2$ is not in $L$. This contradicts the existence of such a function $F$. \\
We still have to generalize the above in two aspects: for an even smaller probability of estimating the volume, and the possibility that the algorithm is non-deterministic. Upon inspection of the proof above, we notice that it can be easily extended in the following way: instead of taking just two families of random bodies, $K_1$ and $K_2$, one may take $d$ different families, $d>2$, which are all indistinguishable by the algorithm, and have different volumes. The
proof can be stretched as far as $d=e^{n^{\frac{\eps}{2}}}$.  To deal with non-determinsitic algorithms, we will use Yao's lemma
(See \cite{RV}, Lemma 11). Let us generate an index $i$ uniformly distributed in $\{1,...,d\}$, then
a body $K$ from the family $K_i$, and then a sequence of uniformly distributed random points on $K$. Following the lines of the above proof,
we see that every deterministic algorithm, given this sequence, will be incorrect in estimating the volume of $K$
with probability (at least) $= 1 - \frac{1}{d} + \mathbf{SE}(n)$.  It follows from Yao's lemma that every non-deterministic algorithm will be incorrect with the same probability for at least one of the families $K_i$. This finishes the theorem.
\qed \\ \\

\bigskip {\noindent School of Mathematical Sciences, Tel-Aviv University, Tel-Aviv
69978, Israel \\  {\it e-mail address:}
\verb"roneneldan@gmail.com" }
\end{document}